\documentclass[12pt,article]{amsart}
\usepackage{mathrsfs}
\numberwithin{equation}{section}
\newtheorem{theorem}{Theorem}[section]
\newtheorem{corollary}{Corollary}
\newtheorem{lemma}[theorem]{Lemma}

 \usepackage{color}

\AtBeginDocument{{\noindent\small
\emph{to appear in Journal of Contemporary Mathematical Analysis}
}
\vspace{8mm}}

\begin{document}
\title[The exact entire solutions of certain type of nonlinear]
{The exact entire solutions of certain type of nonlinear difference equations}
\author[F. L\"{u}, C.P. Li and J. F. Xu]{Feng L\"{u}, Cuiping Li and Junfeng Xu}
\address{Feng L\"{u}\newline
College of Science\\
China University of Petroleum\\
Qingdao, Shandong, 266580, P.R. China.}
\email{lvfeng18@gmail.com}

\address{Cuiping Li\newline
College of Science\\
China University of Petroleum\\
Qingdao, Shandong, 266580, P.R. China.} \email{upclicuiping@163.com}

\address{Junfeng Xu (corresponding author) \newline
Department of Mathematics\\
Wuyi University \\
Jiangmen, Guangdong 529020, P.R.China}
\email{xujunf@gmail.com}


\begin{abstract}
In this paper, we consider the entire solutions of nonlinear difference equation
$$
 f^3+q(z)\Delta f=p_1 e^{\alpha_1 z}+ p_2 e^{\alpha_2 z}
$$
where $q$ is a polynomial, and $p_1, p_2, \alpha_1, \alpha_2$ are nonzero constants with $\alpha_1\neq
\alpha_2$. It is showed that if $f$ is a non-constant entire solution of $\rho_2(f)<1$ to the above
equation, then
$
f(z)=e_1e^{\frac{\alpha_1 z}{3}}+e_2e^{\frac{\alpha_2 z}{3}},
$
where $e_1$ and $e_2$ are two constants. Meanwhile, we give an affirmative answer to
the conjecture posed by  Zhang et al in \cite{Zhang}.
\end{abstract}

\maketitle

\section{\textbf{Introduction and Main results}}
In this paper, by meromorphic functions we will always mean meromorphic
functions in the complex plane. In order to prove the main results, we will employ Nevanlinna theory. Before to proceed, we spare the reader for a moment and assume his/her familiarity with
the basics of Nevanlinna's theory of meromorphic functions in
$\mathbb{C}$ such as the {\it first} and {\it second} fundamental theorems,
and the usual notations such as the {\it characteristic function}
$T(r,f)$, the {\it proximity function} $m(r,f)$ and the {\it
counting function} $N(r,f)$. $S(r,f)$ denotes any quantity
satisfying $S(r,f)=o\left(T(r,f)\right)$ as $r\to\infty$, except
possibly on a set of finite logarithmic measure(see e.g., \cite{BL,Yang,YY}). We also need the following definition.

\noindent \textbf{Definition 1.} The order $\rho(f)$, hyper-order $\rho_{2}(f)$ of the meromorphic function $f(z)$ are defined as follows:
$$\rho (f)=\limsup_{r \rightarrow\infty} \frac{\log T(r,f)}{\log r},~~~~~\rho_{2}(f)=\limsup_{r \rightarrow\infty} \frac{\log\log T(r,f)}{\log r}.$$

Characterizing complex analytic solutions of differential equations has a topic of a long
history (see e.g., the monograph \cite{Hille}). It seems to us that Yang firstly started to study the the existence and uniqueness of finite order entire solutions of nonlinear differential equation of the form
$$
L(f)(z)-p(z)f^n(z)=h(z),~n\geq 3,
$$
where $L(f)$ is a linear differential
polynomial in $f$ with polynomial coefficients, $p$ is a non-vanishing
polynomial and $h$ is an entire function. Recently, the difference analogues to Nevanlinna theory was established
by Halburd and Korhonen \cite{Hal1,HKT}, Chiang and Feng \cite{Chiang}, independently.  With the help of this tool, many scholars have studied the solvability
and existence of meromorphic solutions of some non-linear difference equations (see e.g., \cite{HL,HZ,LY,Latreuch,L,LS,Liu,LLX,LvHan,YL}).

In 2010, Yang and Laine \cite{YL} considered the following difference equation.

{\bf Theorem A}. {\it A non-linear difference equation
\begin{equation*}
 f^3(z)+q(z) f(z+1)=c \sin bz=c \frac{e^{biz}-e^{-biz}}{2i},
\end{equation*}
where $q(z)$ is a non-constant polynomial and $b,c \in \mathbb{C}$ are nonzero constants, does not
admit entire solutions of finite order. If $q(z)=q$ is a nonzero constant, then the above equation possesses
three distinct entire solutions of finite order, provided that $b=3n\pi$ and $q^3=(-1)^{n+1} c^2 27/4$ for a
nonzero integer $n$.}

The follow-up research on this aspect was done by Liu and L\"{u} et al. In \cite{Liu}, they considered the following more general difference equation
\begin{equation}\label{1.01}
 f^n(z)+q(z)\Delta f(z)=p_1 e^{\alpha_1 z}+ p_2 e^{\alpha_2 z},
\end{equation}
where $n$ is a positive integer, $\Delta f(z)=f(z+1)-f(z)$, $q(z)$ is a polynomial, and $p_1, p_2,
\alpha_1, \alpha_2$ are nonzero constants with $\alpha_1\neq \alpha_2$. More specifically, Liu and L\"{u} et al.  proved the following.

{\bf Theorem B}. {\it Let $n\geq 4$ be an integer, $q$ be a polynomial, and $p_1, p_2,
\alpha_1, \alpha_2$ be nonzero constants such that $\alpha_1\neq
\alpha_2$. If there exists some entire solution $f$ of finite order to (\ref{1.01}),
then $q(z)$ is a constant, and one of the following relations holds:\\
\emph{(1)}. $f(z)=c_1e^{\frac{\alpha_1 z}{n}}$, and
$c_1(\exp\frac{\alpha_1}{n}-1)q=p_2, ~\alpha_1=n\alpha_2$,\\
\emph{(2)}. $f(z)=c_2e^{\frac{\alpha_2 z}{n}}$, and
$c_2(\exp\frac{\alpha_2}{n}-1)q=p_1, ~\alpha_2=n\alpha_1$, where
 $c_1,~c_2$ are constants satisfying
$c_1^3=p_1,~c_2^3=p_2$.}

The study for the case $n=3$ was due to Zhang et al. \cite{Zhang}, who obtained the following result.

{\bf Theorem C}. {\it Let $q$ be a polynomial, and $p_1, p_2,
\alpha_1, \alpha_2$ be nonzero constants such that $\alpha_1\neq
\alpha_2$. If $f$ is an entire solution of finite order to the
following equation:
\begin{equation}\label{1.0}
 f^3+q(z)\Delta f=p_1 e^{\alpha_1 z}+ p_2 e^{\alpha_2 z},
\end{equation}
then $q(z)$ is a constant, and one of the following relations holds:

\emph{(1)}. $T(r,f)=N_{1)}(r,\frac{1}{f})+S(r,f)$,

\emph{(2)}. $f(z)=c_1e^{\frac{\alpha_1 z}{3}}$, and
$c_1(\exp\frac{\alpha_1}{3}-1)q=p_2, ~\alpha_1=3\alpha_2$,

\emph{(3)}. $f(z)=c_2e^{\frac{\alpha_2 z}{3}}$, and
$c_2(\exp\frac{\alpha_2}{3}-1)q=p_1, ~\alpha_2=3\alpha_1$, where
$N_{1)}(r,\frac{1}{f})$ denotes the counting function corresponding
to simple zeros of $f$ , and $c_1,~c_2$ are constants satisfying
$c_1^3=p_1,~c_2^3=p_2$.}

\textbf{Remark 1.} For the cases (2) and (3) in Theorem C, it is easy to see that 0 is a Picard value of $f$ and $N(r,1/f)=0$.
So $T(r,f)\neq N_{1)}(r,\frac{1}{f})+S(r,f)=S(r,f)$. It is natural to ask whether the case (1) occurs or not.
The answer is positive. It is showed by the following example, which can be found in \cite{Zhang}.

\textbf{Example 1.} Consider $f(z)=e^{\pi iz}+e^{-\pi iz}=2 i \sin(\pi iz)$.
Then $f$ is a solution of the following equation:
$$
f^3+\frac{3}{2}\Delta f=e^{3\pi iz }+e^{-3\pi iz }.
$$
Obviously, $T(r,f)=N_{1)}(r,\frac{1}{f})+S(r,f)$. So, the case (1) occurs.

In Theorem C, it seems that the case (1) is unnatural. Meanwhile, Zhang et al. observed that $\alpha_1+\alpha_2=3\pi i+(-3\pi i)=0$ in Example 1. This observation leaded Zhang et al. to pose the following conjecture.

\textbf{Conjecture.} {\it If $\alpha_1\neq \alpha_2$, $\alpha_1+\alpha_2\neq 0$, then the conclusion (1) of Theorem C is
impossible. In fact, any entire solution f of (\ref{1.0}) must have 0 as its Picard exceptional value.}

\textbf{Remark 2.} The conjecture has been studied by many researcher.
In 2017, Latreuch in \cite{Latreuch} has gave an affirmative answer to the conjecture.
However, when the $\alpha_1+\alpha_2\neq 0$ does not hold,
Latreuch did not give the specific form of the meromorphic solution of (\ref{1.0}).
In Example 1, we further observe that $f(z)=e^{\pi iz}+e^{-\pi iz}=2 i \sin(\pi iz)$.  This leaded us to ask whether any entire solution of the equation (\ref{1.0}) always is this form when Case (1) occurs. In the present paper, we focus on the problem and give an affirmative answer by the following theorem.

\begin{theorem}\label{th1}
Let $q$ be a polynomial, and $p_1, p_2,
\alpha_1, \alpha_2$ be nonzero constants such that $\alpha_1\neq
\alpha_2$. If $f$ is an entire solution of $\rho_2(f)<1$ to the
equation (\ref{1.0}), then $q(z)$ is a constant, and one of the following
relations holds:

\emph{(1)}. $f(z)=e_1e^{\frac{\alpha_1 z}{3}}+e_2e^{\frac{\alpha_2 z}{3}}$,
where $e_1$ and $e_2$ are two nonzero constants satisfying
$e_1^3=p_1$, $e_2^3=p_2$ (or $e_1^3=p_2$, $e_2^3=p_1$),
$3e_1e_2-2q=0$, $\alpha_1+\alpha_2=0$ and
$e^{\frac{\alpha_1}{3}}=-1$;

\emph{(2)}. $f(z)=c_1e^{\frac{\alpha_1 z}{3}}$, and
$c_1(\exp\frac{\alpha_1}{3}-1)q=p_2, ~\alpha_1=3\alpha_2$;

\emph{(3)}. $f(z)=c_2e^{\frac{\alpha_2 z}{3}}$, and
$c_2(\exp\frac{\alpha_2}{3}-1)q=p_1, ~\alpha_2=3\alpha_1$.
\end{theorem}


\textbf{Remark 3.} Clearly, Example 1 satisfies Case (1) of Theorem 1.1, where $\alpha_1=3\pi i, \alpha_2=-3\pi i$; $e_1=e_2=1, p_1=p_2=1$; $q=3/2$.
Next we give two examples to show Cases (2) and (3) indeed occur in Theorem \ref{th1}.

\textbf{Example 2.} Consider the function $f(z)=e^{\pi iz}$, which is a nonconstant entire solution of the following equation
$$
f^3(z)-\frac{1}{2}\Delta f(z)=e^{3\pi iz }+e^{\pi iz },
$$
where $\alpha_1=3\pi i=3\alpha_2$, $c_1=1,q=-1/2, p_2=1$. Thus, the case (2) occurs.

\textbf{Example 3.} Consider the function $f(z)=e^{3\pi iz}$, which satisfies the following equation
$$
f^3(z)-\frac{1}{2}\Delta f(z)=e^{3\pi iz }+e^{9\pi iz },
$$
where $\alpha_2=9\pi i=3\alpha_1$, $c_2=1,q=-1/2, p_1=1$. Therefore, the case (3) occurs.

By Theorem \ref{th1}, we get an immediate conclusion as follows.

\begin{corollary}\label{co1}
Let $q$ be a polynomial, and $p_1, p_2,
\alpha_1, \alpha_2$ be nonzero constants such that $\alpha_1\neq
\alpha_2$. If $f$ is a nonconstant entire solution of $\rho_2(f)<1$ to the
equation (\ref{1.0}), then $q(z)$ is a constant, and
$$
f(z)=e_1e^{\frac{\alpha_1 z}{3}}+e_2e^{\frac{\alpha_2 z}{3}},
$$
where $e_1$ and $e_2$ are two constants.
\end{corollary}

At the end, we turn attention to the question: What will happen if we replace the function $f^3$ by $f^2$ in the equation (\ref{1.0}).
After studying this question, we derive some similar results to Theorem C as follows.

\begin{theorem}\label{th2}
Let $q$ be a polynomial, and $p_1, p_2,
\alpha_1, \alpha_2$ be nonzero constants such that $\alpha_1\neq
\alpha_2$. If $f$ is an entire solution of $\rho_2(f)<1$ to the
following equation
\begin{equation}\label{1.1}
 f^2+q(z)\Delta f=p_1 e^{\alpha_1 z}+ p_2 e^{\alpha_2 z},
\end{equation}
and  satisfying $N(r,\frac{1}{f})=S(r,f)$, then $q(z)$ is a constant, and one of the following relations holds:

\emph{(1)}. $f(z)=c_1e^{\frac{\alpha_1 z}{2}}$, and
$c_1(\exp\frac{\alpha_1}{2}-1)q=p_2, ~\alpha_1=2\alpha_2$,

\emph{(2)}. $f(z)=c_2e^{\frac{\alpha_2 z}{2}}$, and
$c_2(\exp\frac{\alpha_2}{2}-1)q=p_1, ~\alpha_2=2\alpha_1$, where $c_1,~c_2$ are constants satisfying
$c_1^2=p_1,~c_2^2=p_2$.
\end{theorem}

We below offer an example to show that the condition $N(r,\frac{1}{f})=S(r,f)$ is necessary in Theorem \ref{th2}.

\textbf{Example 4.} Consider the function $f(z)=-2-\sqrt{2}e^{\pi iz}+\sqrt{2}e^{-\pi iz}$, which satisfies the equation
$$
f^2(z)-2\Delta f(z)=2e^{2\pi iz }+2e^{-2\pi iz }.
$$
A calculation yields that $T(r,f)=2r(1+o(1))$ and $N(r,1/f)=2r(1+o(1))$. Clearly, $N(r,\frac{1}{f})\neq S(r,f)$ and $f$ does not satisfy any  conclusion of Theorem \ref{th2}.

\section{\textbf{Some Lemmas}}

Before to the proofs of main theorems, we
firstly give the following result, whcih is a version of the difference analogue of the logarithmic derivative lemma.

\begin{lemma}[\cite{HKT}]\label{le1}
Let $f(z)$ be a meormorphic function of $\rho_2(f)<1$, and let
$c\in \mathbb{C}\backslash \{0\}$. Then
$$
    m(r, \frac{f(z+c)}{f(z)})=o(\frac{T(r,f)}{r^{1-\rho_2(f)-\varepsilon}}),
$$
outside of an exceptional set of finite logarithmic
measure.
\end{lemma}

In addition, by applying Lemma 2.1 and the same argument as in \cite[Theorem 2.3]{LY}, we get
the following lemma, which is a version of the difference analogue of the Clunie lemma. The details are omitted here.

\begin{lemma}\label{le2}
Let $f$ be a transcendental meromorphic solution of  $\rho_2(f)<1$
to the difference equation
$$
H(z,f)P(z,f)=Q(z,f),
$$
where $H(z,f),~P(z,f),~Q(z,f)$ are difference polynomials in $f$ such
that the total degree of $H(z,f)$ in $f$ and its shifts is $n$, and
that the corresponding total degree of $Q(z,f)$ is $\leq n$. If
$H(z,f)$ contains just one term of maximal total degree, then for
any $\varepsilon>0$,
$$
m(r,P(z,f))=S(r,f),
$$
possibly outside of an exceptional set of finite logarithmic
measure.
\end{lemma}

\section{\textbf{Proof of Theorem \ref{th1}}}

Suppose that $f$ is an entire solution of $\rho_2(f)<1$ to Eq (\ref{1.0}). Obviously, $f$ is a transcendental function. By
differentiating both sides of (\ref{1.0}), one has
\begin{equation}\label{2.1}
3f^2f'+(q(z)\Delta f)'=\alpha_1 p_1 e^{\alpha_1 z}+ \alpha_2p_2
e^{\alpha_2 z}.
\end{equation}
Combining (\ref{1.0}) and (\ref{2.1}) yields
\begin{equation}\label{2.3}
\alpha_2f^3+\alpha_2 q \Delta f- 3f^2f'-(q(z)\Delta
f)'=(\alpha_2-\alpha_1) p_1 e^{\alpha_1 z}.
\end{equation}
By differentiating (\ref{2.3}), we derive that
\begin{equation}\label{2.4}
3\alpha_2f^2f'+\alpha_2 (q \Delta f)'- 6f(f')^2-3f^2f''-(q(z)\Delta
f)''=\alpha_1(\alpha_2-\alpha_1) p_1 e^{\alpha_1 z}.
\end{equation}
It follows from (\ref{2.3}) and (\ref{2.4}) that
\begin{equation}\label{2.5}
f\varphi=T(z,f),
\end{equation}
where
\begin{equation}\label{2.6}
\varphi=\alpha_1 \alpha_2 f^2-3(\alpha_1+\alpha_2)ff'+6(f')^2+3ff'',
\end{equation}
$$
T(z,f)=-\alpha_1 \alpha_2 q \Delta f+(\alpha_1+\alpha_2)(q \Delta
f)'-(q \Delta f)''.
$$
Note that $T(z,f)$ is a differential-difference polynomial in $f$ of
degree 1. Then by applying Lemma \ref{le2} to the equation (\ref{2.5}), one has
$m(r,\varphi)=S(r,f)$. Further, $T(r,\varphi)=m(r,\varphi)=S(r,f)$, since $\varphi$ is an entire function. It means that $\varphi$ is a
small function of $f$.

 Suppose that $\varphi\equiv 0$. Then $\alpha_1 \alpha_2 f^2-3(\alpha_1+\alpha_2)ff'+6(f')^2+3ff''\equiv0$. Rewrite it as  $\frac{f''}{f}=(\frac{f'}{f})'+(\frac{f'}{f})^2$,  which yields a Riccati equation
$$ t'+3t^2-(\alpha_1+\alpha_2)t+\alpha_1\alpha_2/3=0,$$
where $t=\frac{f'}{f}$. Clearly, the equation  has two constant solutions $t_1=\alpha_1/3$, $t_2=\alpha_2/3$. We assume $t\not\equiv t_1, t_2$. Then we have
\begin{equation*}
  \frac{1}{t_1-t_2}(\frac{t'}{t-t_1}-\frac{t'}{t-t_2})=-3.
\end{equation*}
Integrating the above equation yields
\begin{equation*}
 \ln \frac{t-t_1}{t-t_2}=3(t_2-t_1)z+C,
\end{equation*}
where $C$ is a constant. Therefore,
\begin{equation*}
 \frac{t-t_1}{t-t_2}=e^{3(t_2-t_1)z+C}.
\end{equation*}
This immediately yields
\begin{equation*}
 t=t_2+\frac{t_2-t_1}{e^{3(t_2-t_1)z+C}-1}=\frac{f'}{f},
\end{equation*}

Note that the zeros of $e^{3(t_2-t_1)z+C}-1$ are the zeros of $f$. If $z_0$ is a zero of $f$ with the multiplicity $k$,
then
$$ k=\text{Res} [\frac{f'}{f}, z_0] =\text{Res} [t_2+\frac{t_2-t_1}{e^{3(t_2-t_1)z+C}-1}, z_0] =\frac{1}{3},$$
which is a contradiction. Thus, either $t\equiv t_1=\alpha_1/3$ or $t\equiv t_2=\alpha_2/3$.

If $t\equiv t_1=\alpha_1/3$, then $f(z)=c_1 e^{{\frac{\alpha_1}{3}z}}$. Substituting the form $f(z)=c_1 e^{{\frac{\alpha_1}{3}z}}$ into the equation (\ref{1.0}), we obtain
$$
c_1^3 e^{\alpha_1 z}+c_1 q(z)e^{\frac{\alpha_1}{3} z}(e^{\frac{\alpha_1}{3}}-1)= p_1 e^{\alpha_1 z}+ p_2 e^{\alpha_2 z},$$
which implies that $c_1^3=p_1$, $c_1 q(e^{\frac{\alpha_1}{3}}-1)=p_2$ and $\alpha_1=3\alpha_2$.

Similarly as above, if $t\equiv t_2=\alpha_2/3$, then we can derive that $f(z)=c_2 e^{{\frac{\alpha_2}{3}z}}$ satisfying $c_2^3=p_2$, $c_2 q(e^{\frac{\alpha_2}{3}}-1)=p_1$ and $\alpha_2=3\alpha_1$.\\


In the following, based on the idea in \cite[Theorem 1.1]{L},
we will consider the case $\varphi\not\equiv 0$. By Theorem C, one has
\begin{equation}\label{2.7}
T(r,f)=N_{1)}(r,\frac{1}{f})+S(r,f).
\end{equation}
Differentiating (\ref{2.6}) yields
\begin{equation}\label{2.8}
\varphi'=\alpha_1 \alpha_2
2ff'-3(\alpha_1+\alpha_2)(ff''+(f')^2)+12f'f''+3ff'''+3f'f''.
\end{equation}

From (\ref{2.6}) and (\ref{2.8}), we can obtain that
\begin{equation}\label{2.9}
f[A_0f+A_1f'+A_2f''+A_3f''']=f'[B_1f'+B_2f''],
\end{equation}
where
$$
\begin{aligned}
&A_0=
\alpha_1\alpha_2\varphi',~~A_1=-3\varphi'(\alpha_1+\alpha_2)-2\varphi\alpha_1\alpha_2,\\
&A_2=3\varphi'+3\varphi(\alpha_1+\alpha_2),~~A_3=-3\varphi,\\
&B_1=-3\varphi(\alpha_1+\alpha_2)-6\varphi',~~B_2=15\varphi.
\end{aligned}
$$
Obviously, all $A_i$ $(i=0,1,2,3)$, $B_j$ $(j=1,2)$ are small
functions of $f$.

Suppose that $z_0$ is a zero of $f$, not a zero of $\varphi$. It
follows from (\ref{2.6}) that $6(f')^2(z_0)=\varphi(z_0)\neq 0$,
which implies that $z_0$ is a simple zero of $f$. Then by
(\ref{2.9}), we have
$$
B_1(z_0)f'(z_0)+B_2(z_0)f''(z_0)=0.
$$
Set
\begin{equation}\label{2.10}
A=\frac{B_1f'+B_2f''}{f}.
\end{equation}
We claim that $A$ is an entire function. Clearly, all the simple
zeros of $f$ are not poles of $f$. Suppose that $b_0$ is a multiple
zero of $f$. By (\ref{2.6}), we get $b_0$ is also a multiple zero of
$\varphi$. So, $b_0$ is a zero of $B_1$ and a multiple zero of
$B_2$. Note that $b_0$ is a pole of $\frac{f'}{f}$ and
$\frac{f''}{f}$ with multiplicity one and two, respectively. Thus,
$b_0$ is not a pole of $B_1\frac{f'}{f}$ and $B_2\frac{f''}{f}$,
which implies that $b_0$ is not a pole of $A$. Thus, $A$ is an
entire function. The claim is proved. Furthermore,
$$
T(r, A)=m(r,\frac{B_1f'+B_2f''}{f})=S(r,f).
$$
Hence $A$ is a small function of $f$. We consider two cases below.

Case 1. $A=0$.

Then, $B_1f'+B_2f''=0$. Rewrite it as
$$
\frac{f''}{f'}=-\frac{B_1}{B_2}=\frac{1}{5}(\alpha_1+\alpha_2)+\frac{2}{5}\frac{\varphi'}{\varphi}.
$$
By integrating the above equation, we have
$$
f'(z)=\beta e^{\frac{1}{5}(\alpha_1+\alpha_2)z},
$$
where $\beta$ is a small function of $f$. Obviously,
$\alpha_1+\alpha_2\neq 0$. Otherwise, $T(r,f')=T(r,\beta)=S(r,f)$, a
contradiction. We below consider two subcases.

Subcase 1.1. $\varphi'=0$.

The equation $B_1f'+B_2f''=0$ yields
$$
\frac{f''}{f'}=-\frac{B_1}{B_2}=\frac{1}{5}(\alpha_1+\alpha_2).
$$
By integrating the above equation, we derive that
$f'(z)=H_1e^{\frac{1}{5}(\alpha_1+\alpha_2)z}$, where $H_1$ is a
nonzero constant.

Integrating the function $f'$ yields
$$
f(z)=k_1e^{\frac{1}{5}(\alpha_1+\alpha_2)z}+k_2,
$$
where $k_1(\neq 0),~~k_2$ are two constants. Obviously, $k_2\neq0$.
Otherwise, $f$ has no zeros, which contradicts with (\ref{2.7}).
Substitute the form of $f$ into the equation (\ref{1.0}) yields
$$
\begin{aligned}
&a_3e^{\frac{3}{5}(\alpha_1+\alpha_2)z}+a_2e^{\frac{2}{5}(\alpha_1+\alpha_2)z}\\
&+a_1e^{\frac{1}{5}(\alpha_1+\alpha_2)z}+k_2^3=p_1 e^{\alpha_1 z}+
p_2 e^{\alpha_2 z},
\end{aligned}
$$
where $a_1,a_2,a_3$ are small functions of $f$. Then, the above
equation yields that $k_2=0$, a contradiction. Hence Subcase 1.1 can
not occur.

Subcase 1.2.  $\varphi'\neq0$.

By differentiating $f'$ one and two times respectively, we have
$$
f''=H_2e^{\frac{1}{5}(\alpha_1+\alpha_2)z},~~f'''=H_3e^{\frac{1}{5}(\alpha_1+\alpha_2)z},
$$
where $H_2$ and $H_3$ are two small functions of $f$. The equation
(\ref{2.9}) implies that
$$
A_0f+A_1f'+A_2f''+A_3f'''=0.
$$
Furthermore,
$$
f=-\frac{A_1f'+A_2f''+A_3f'''}{A_0}=H_0
e^{\frac{1}{5}(\alpha_1+\alpha_2)z},
$$
where $H_0$ is a small function of $f$. So,
$$
N(r,\frac{1}{f})=N(r,\frac{1}{H_0})\leq T(r,H_0)=S(r,f),
$$
which contradicts with (\ref{2.7}). Thus, Subcase 1.2 can not
occur.

Case 2. $A\neq0$.

By (\ref{2.9}) and (\ref{2.10}),  one has
$$
\frac{A_0f+A_1f'+A_2f''+A_3f'''}{f'}=A,
$$
which yields that
\begin{equation}\label{2.11}
A_0f+(A_1-A)f'+A_2f''+A_3f'''=0.
\end{equation}
Rewrite (\ref{2.10}) as
$$
Af-B_1f'-B_2f''=0.
$$
Differentiating the above equation as
\begin{equation}\label{2.12}
A'f+(A-B_1')f'-(B_1+B_2')f''-B_2f'''=0.
\end{equation}
Combining (\ref{2.11}) and (\ref{2.12}) yields
\begin{equation}\label{2.A}
C_0f+C_1f'+C_2f''=0,
\end{equation}
where
$$
\begin{aligned}
C_0= A_0B_2+A'A_3,~~
C_1=(A_1-A)B_2+A_3(A-B_1'),~~C_2=A_2B_2-A_3(B_1+B_2').
\end{aligned}
$$
Obviously, $C_i$ $(i=0,1,2)$ are small functions of $f$.

We consider two subcases again.

Subcase 2.1. $C_2=0$.

It follows that $C_0=C_1=0$. Otherwise, without loss of generality, suppose
that $C_0\neq 0$. By (\ref{2.A}), we have that $C_1\neq0$. Assume
that $\omega_0$ is a simple zero of $f$. Then $\omega_0$ is a zero
of $C_1$. Furthermore,
$$
T(r,f)=N_{1)}(r,\frac{1}{f})+S(r,f)\leq
N(r,\frac{1}{C_1})+S(r,f)\leq T(r, C_1)+S(r,f)=S(r,f),
$$
a contradiction. Thus, $C_0=C_1=0$.

The fact $C_2=0$ leads to
\begin{equation}\label{2.13}
2\varphi'+\varphi(\alpha_1+\alpha_2)=0.
\end{equation}
If $\alpha_1+\alpha_2\neq 0$, then
$\varphi=H_4e^{-\frac{\alpha_1+\alpha_2}{2}z}$, where $H_4$ is a
nonzero constant. Therefore, we have
$$
\begin{aligned}
&m(r,\varphi)=\frac{|\frac{\alpha_1+\alpha_2}{2}|}{\pi}r(1+o(1)),\\
&m(r, e^{\alpha_1 z})=\frac{|\alpha_1|}{\pi}r(1+o(1)), \\
&m(r, e^{\alpha_2 z})=\frac{|\alpha_2|}{\pi}r(1+o(1)).
\end{aligned}
$$
Note that $\varphi$ is a small function of $f$. So $e^{\alpha_1
z},~e^{\alpha_2 z}$ are also two small functions of $f$. Rewrite
(\ref{1.0}) as
$$
f^3=-q(z)\Delta f+p_1 e^{\alpha_1 z}+ p_2 e^{\alpha_2 z}.
$$
Therefore,
$$
\begin{aligned}
3T(r,f)&=T(r,f^3)=T(r,-q(z)\Delta f+p_1 e^{\alpha_1 z}+ p_2
e^{\alpha_2 z})\\
&\leq T(r, \Delta f)+S(r,f)\leq T(r,f)+S(r,f),
\end{aligned}
$$
a contradiction.

Hence $\alpha_1+\alpha_2=0$. Then, (\ref{2.13}) reduces to
$\varphi'=0$. It implies that $\varphi$ is a constant and
$A_0=\varphi'\alpha_1 \alpha_2=0$. Together with $C_0=0$, it is
easy to deduce that $A'=0$ and $A$ is also a constant. Therefore, $B_1$ and
$B_2$ become two constants. Then the following equation reduces to a constant
coefficient homogeneous linear differential equation
$$
Af-B_1f'-B_2f''=0.
$$

Suppose that the characteristic equation $B_2\lambda^2+B_1\lambda-A=0$ has two distinct roots
$\lambda _1, ~\lambda_2$. Clearly, $\lambda _1, ~\lambda_2$ are
nonzero constants. Then, by solving the above equation, one derives
\begin{equation}\label{2.14}
f(z)=e_1 e^{\lambda_1 z}+ e_2 e^{\lambda_2 z}.
\end{equation}
Clearly, $e_1e_2\neq 0$. Otherwise $f$ has no zeros, a
contradiction. Substitute the form $f$ into (\ref{1.0}), we have
\begin{equation}\label{B}
\begin{aligned}
e_1^3 e^{3\lambda_1 z }+e_2^3 e^{3\lambda_2 z
}+&3e_1^2e_2e^{(2\lambda_1+\lambda_2)z}+3e_1e_2^2e^{(\lambda_1+2\lambda_2)z}\\
&+qe_1(e^{\lambda_1}-1) e^{\lambda_1 z}+qe_2(e^{\lambda_2}-1)e^{\lambda_2
z}=p_1e^{\alpha_1 z}+p_2e^{\alpha_2 z}.
\end{aligned}
\end{equation}
Suppose that $\lambda_1+\lambda_2\neq 0$. Observe that $\lambda_1\neq
\lambda_2$. So $3\lambda_1,~3\lambda_2,~2\lambda_1+\lambda_2,~\lambda_1+2\lambda_2$
are distinct from each other. Furthermore, by (\ref{B}) and Borel's Theorem,
we easily get the following two sets are identity
$$
\{3\lambda_1,~3\lambda_2,~2\lambda_1+\lambda_2,~\lambda_1+2\lambda_2\}=\{\lambda_1,~\lambda_2,~\alpha_1,~\alpha_2\},
$$
which implies that $3\lambda_2=\lambda_1$ and
$3\lambda_1=\lambda_2$. It is impossible. Thus,
$\lambda_1+\lambda_2=0$. Rewrite (\ref{B}) as
$$
e_1^3 e^{3\lambda_1 z }+e_2^3 e^{3\lambda_2 z }+q_1 e^{\lambda_1
z}+q_2e^{\lambda_2 z}=p_1e^{\alpha_1 z}+p_2e^{\alpha_2 z},
$$
where $q_1=3e_1^2e_2+q e_1(e^{\lambda_1}-1)$, $q_2=3e_2^2e_1+q
e_2(e^{\lambda_2}-1)$ are two polynomials. Then, it follows from the
above equation that $q_1=q_2=0$. Meanwhile, one has
$$
3\lambda_1= \alpha_1,~~3\lambda_2=\alpha_2
$$
or
$$ 3\lambda_1= \alpha_2,~~3\lambda_2=\alpha_1.
$$
Furthermore, we obtain that $e_1^3=p_1$ and $e_2^3=p_2$ (or
$e_1^3=p_2$ and $e_2^3=p_1$). Note that
$$
q_1=3e_1^2e_2+qe_1(e^{\lambda_1}-1)=0,~~q_2=3e_2^2e_1+qe_2(e^{\lambda_2}-1)=0.
$$
By the above two equation, $\lambda_1+\lambda_2=0$ and a calculation, we deduce that $e^{\lambda_1}=-1$ and $q$ reduces to a constant
satisfying $3e_1e_2-2q=0$.\\

Now, we suppose that $B_2\lambda^2+B_1\lambda-A=0$ has a multiple root,
say $\lambda_3$. Then, $f(z)=(e_3+e_4 z)e^{\lambda_3 z}$. Therefore, $f$
just has one zero, a contradiction.

Subcase 2.2. $C_2\neq 0$.

Combining (\ref{2.10}) and (\ref{2.A}) yields
$$
(B_2C_0+AC_2)f+(C_1B_2-B_1C_2)f'=0.
$$
Suppose that $C_1B_2-B_1C_2\neq 0$. It folllows $B_2C_0+AC_2\neq0$. Assume
that $\sigma_0$ is a simple zero of $f$. By the above equation, one has $\sigma_0$ is also a zero of $C_1B_2-B_1C_2$. Then,
$$
\begin{aligned}
T(r,f)&=N_{1)}(r,\frac{1}{f})+S(r,f)\leq N(r,
\frac{1}{C_1B_2-B_1C_2})+S(r,f)\\
&\leq T(r, C_1B_2-B_1C_2)+S(r,f)=S(r,f),
\end{aligned}
$$ a contradiction. The above discussion forces that $C_1B_2-B_1C_2=0$ and $B_2C_0+AC_2=0$.
By the definitions of $C_1,~C_2,~B_1,~B_2$,  a calculation leads to
\begin{equation}\label{2.15}
8A\varphi'-5\varphi A' =-[4\varphi
A(\alpha_1+\alpha_2)+25\alpha_1\alpha_2\varphi\varphi']
\end{equation}
and
\begin{equation}\label{2.16}
15\varphi
A=[6(\alpha_1+\alpha_2)^2-25\alpha_1\alpha_2]\varphi^2-21(\alpha_1+\alpha_2)\varphi\varphi'+24(\varphi')^2-15\varphi\varphi''.
\end{equation}
Suppose that $\delta_0$ is a zero of $\varphi$ with multiplicity
$s$. The equation (\ref{2.16}) implies $s\geq 2$. Furthermore, $\delta_0$ is a zero of $\varphi^2$ and
$\varphi\varphi'$ with multiplicity $2s$ and $2s-1$, respectively.
Suppose that the Laurent expansions of $\varphi$  at $\delta_0$ is as follows
$$\varphi(z)=\mu_s(z-\delta_0)^s+\mu_{s+1}(z-\delta_0)^{s+1}+\cdots,
$$
where $\mu_s(\neq 0),~~\mu_{s+1}$ are constants. Then, a calculation
yields
$$
24(\varphi')^2-15\varphi\varphi''=[24(\mu_s)^2s^2-15(\mu_s)^2s(s-1)](z-\delta_0)^{2s-2}+\theta_{2s-1}(z-\delta_0)^{2s-1}+\cdots,
$$
where $\theta_{2s-1}$ is a constant. Obviously,
$$
24\mu_s^2s^2-15\mu_s^2s(s-1)=\mu_s^2s[9s+15]\neq 0,
$$
which implies that $\delta_0$ is a zero of
$24(\varphi')^2-15\varphi\varphi''$ with multiplicity $2s-2$.
Suppose that $\delta_0$ is a zero of $A$ with multiplicity $l$.
Then, comparing the multiplicity of both side of equation
(\ref{2.16}) at point $\delta_0$, we have $s+l=2s-2$. So, $s=l+2$.

Assume that $l=0$. Then, $s=2$ and $A(\delta_0)\neq 0$. Rewrite
(\ref{2.15}) as
\begin{equation}\label{2.17}
8A\varphi'=5\varphi A' -[4\varphi
A(\alpha_1+\alpha_2)+25\alpha_1\alpha_2\varphi\varphi'].
\end{equation}
Clearly, $\delta_0$ is a simple zero of $A\varphi'$. However,
$\delta_0$ is a multiple zero of $5\varphi A' -[4\varphi
A(\alpha_1+\alpha_2)+25\alpha_1\alpha_2\varphi\varphi']$, a
contradiction. Therefore, $l\geq 1$.

Furthermore, $\delta_0$ is a zero of $4\varphi
A(\alpha_1+\alpha_2)+25\alpha_1\alpha_2\varphi \varphi'$ with
multiplicity $2l+2$. Suppose that the Laurent expansions of $A$  at $\delta_0$ is
$$A(z)=\nu_l(z-\delta_0)^l+\nu_{s+1}(z-\delta_0)^{l+1}+\cdots,
$$
Then,
$$
8A\varphi'-5\varphi A'=\nu_l
\mu_{l+2}[8(l+2)-5l](z-\delta_0)^{2l+1}+\xi_{2l+2}(z-\delta_0)^{2l+2}+\cdots,
$$
where $\xi_{2l+2}$ is a constant. Then, $\delta_0$ is a zero of
$8A\varphi'-5\varphi A'$ with multiplicity $2l+1$, since $\nu_l
\mu_{l+2}[8(l+2)-5l]\neq 0$. So, the point $\delta_0$ is a zero of the left side function of (\ref{2.15}) with multiplicity $2l+1$. On the other hand, $\delta_0$ is a zero of the right side function of (\ref{2.15}) with multiplicity at least  $2l+2$, which is impossible. Therefore, $\varphi$ has no zeros.

If $\varphi$ is not a constant, then, we can assume that $\varphi=
\phi e^{\omega (z)}$, where $\phi$ is a constant and $\omega(\neq
0)$ is an entire function. Then, the same argument as in Subcase 2.1 yields that $e^{\alpha_1 z}$ and $e^{\alpha_2 z}$ are two small functions
of $f$. Furthermore, we can derive a contradiction. Thus, $\varphi$ is a constant. Plus (\ref{2.16}), one has that
$A$ is also a constant. Furthermore, it follows from (\ref{2.15}) that
$\alpha_1+\alpha_2=0$. Similarly as the above discussion, we can deduce the desired result.

Thus, we finish the proof of Theorem \ref{th1}.

\section{\textbf{Proof of Theorem \ref{th2}}}

Suppose that $f$ is an entire solution of $\rho_2(f)<1$ to the equation
(\ref{1.1}). Obviously, $f$ is a transcendental function. By
differentiating both sides of (\ref{1.1}), one has
\begin{equation}\label{3.1}
2f f'+(q(z)\Delta f)'=\alpha_1 p_1 e^{\alpha_1 z}+ \alpha_2p_2
e^{\alpha_2 z}.
\end{equation}
Combining (\ref{1.1}) and (\ref{3.1}) yields
\begin{equation}\label{3.2}
\alpha_2f^2+\alpha_2 q \Delta f- 2ff'-(q(z)\Delta
f)'=(\alpha_2-\alpha_1) p_1 e^{\alpha_1 z}.
\end{equation}
By differentiating (\ref{3.2}), we derive that
\begin{equation}\label{3.3}
2\alpha_2ff'+\alpha_2 (q \Delta f)'- 2(f')^2-2ff''-(q(z)\Delta
f)''=\alpha_1(\alpha_2-\alpha_1) p_1 e^{\alpha_1 z}.
\end{equation}
It follows from (\ref{3.2}) and (\ref{3.3}) that
\begin{equation}\label{3.4}
\varphi_1=T_1(z,f),
\end{equation}
where
\begin{equation}\label{3.5}
\varphi_1=\alpha_1 \alpha_2 f^2-2(\alpha_1+\alpha_2)ff'+2ff''+2(f')^2,
\end{equation}
$$
T_1(z,f)=-\alpha_1 \alpha_2 q \Delta f+(\alpha_1+\alpha_2)(q \Delta
f)'-(q \Delta f)''.
$$
If $\varphi_1\not\equiv 0$, then
$$
\frac{1}{f^2}=\frac{1}{\varphi_1}(\alpha_1 \alpha_2 -2(\alpha_1+\alpha_2)\frac{f'}{f}+2\frac{f''}{f}+2(\frac{f'}{f})^2).
$$
By (\ref{3.4})-(\ref{3.5}), and Lemma 2.1, we have
\begin{equation}\label{3.6}
m(r,\frac{\varphi_1}{f})=m(r,\frac{T_1}{f})=S(r,f) \indent \text{and}\indent m(r,\frac{\varphi_1}{f^2})=S(r,f).
\end{equation}
Combining $N(r,\frac{1}{f})=S(r,f)$ and (\ref{3.6}), we obtain
\begin{eqnarray}
\nonumber   2T(r,f) =2 m(r,\frac{1}{f})+S(r,f)& = & m(r, \frac{\varphi_1}{f^2})+S(r,f) \\
\nonumber   &\leq & m(r, \frac{\varphi_1}{f^2}) + m(r, \frac{1}{\varphi_1})+S(r,f)\\
\nonumber   &\leq & T(r,\varphi_1) +S(r,f)= m(r, \varphi_1) +S(r,f) \\
\nonumber   &=& m(r, \frac{\varphi_1}{f}) + m(r, f)+S(r,f) \\
\nonumber   &=& T(r,f)+S(r,f),
\end{eqnarray}
which implies $T(r,f)=S(r,f)$, a contradiction.

If $\varphi_1\equiv 0$, then by the similar reasoning as in Theorem
1.1 we can obtain the conclusions (1) and (2). Below, we give the details. By $\varphi_1\equiv 0$, one has the differential equation $\alpha_1 \alpha_2 f^2-2(\alpha_1+\alpha_2)ff'+2ff''+2(f')^2=0$. Plus the fact $\frac{f''}{f}=(\frac{f'}{f})'+(\frac{f'}{f})^2$, we can rewrite the above equation to a Riccati equation
$$ t'+2t^2-(\alpha_1+\alpha_2)t+\alpha_1\alpha_2/2=0,$$
where $t=\frac{f'}{f}$. Clearly, the equation has two constant solutions $t_1=\alpha_1/2$, $t_2=\alpha_2/2$.

Suppose the solution $t\not\equiv t_1, t_2$. Then
\begin{equation*}
  \frac{1}{t_1-t_2}(\frac{t'}{t_1-t_2}-\frac{t'}{t_1-t_2})=-2.
\end{equation*}
Integrating the above equation yields
\begin{equation*}
 \ln \frac{t-t_1}{t-t_2}=2(t_2-t_1)z+C,
\end{equation*}
where $C$ is a constant. Therefore,
\begin{equation*}
 \frac{t-t_1}{t-t_2}=e^{2(t_2-t_1)z+C}.
\end{equation*}
This immediately yields
\begin{equation*}
 t=t_2+\frac{t_2-t_1}{e^{2(t_2-t_1)z+C}-1}=\frac{f'}{f}.
\end{equation*}

Note that the zeros of $e^{2(t_2-t_1)z+C}-1$ are the zeros of $f$. If $z_0$ is the zero of $f$ with the multiplicity $k$,
then
$$ k=\text{Res}[\frac{f'}{f}, z_0] =\text{Res} [t_2+\frac{t_2-t_1}{e^{2(t_2-t_1)z+C}-1}, z_0] =\frac{1}{2}.$$
It is a contradiction. Thus, either $t\equiv t_1=\alpha_1/2$ or $t\equiv t_2=\alpha_2/2$.

If $t\equiv t_1=\alpha_1/2$, then $f(z)=c_1 e^{{\frac{\alpha_1}{2}z}}$. Substituting $f(z)=c_1 e^{{\frac{\alpha_1}{2}z}}$ into (\ref{1.1}), we obtain
$$
c_1^2 e^{\alpha_1 z}+c_1 q(z)e^{\frac{\alpha_1}{2} z}(e^{\frac{\alpha_1}{2}}-1)= p_1 e^{\alpha_1 z}+ p_2 e^{\alpha_2 z}. $$
Moreover, we have $c_1^2=p_1$, $c_1 q(e^{\frac{\alpha_1}{2}}-1)=p_2$ and $\alpha_1=2\alpha_2$.

Similarly, if $t\equiv t_2=\alpha_2/2$, then we have $f(z)=c_2 e^{{\frac{\alpha_2}{2}z}}$
satisfying $c_2^2=p_2$, $c_2 q(e^{\frac{\alpha_2}{2}}-1)=p_1$ and
$\alpha_2=2\alpha_1$.

Thus, we finish the proof of Theorem \ref{th2}.

\section*{Acknowledgements}
The research was supported by  Guangdong Basic and Applied Basic Reserch Foundation(No.2018A0303130058),
NNSF of China (Nos. 11601521, 11871379), Funds of Education Department of Guangdong (2016KTSCX145, 2019KZDXM025)
 and the Fundamental Research Fund for Central Universities in China Project(No. 18CX02048A).


\begin{thebibliography}{99}


\bibitem{BL}W. Bergweiler, J.K. Langley, ``Zeros of differences of meromorphic functions",
{\it  Math. Proc. Camb. Phil. Soc.} {\bf142}, 133--147, 2007.

\bibitem{Chiang} Y.M. Chiang and S.J. Feng, ``On the Nevanlinna characteristic of $f(z +\eta)$ and difference equations in
the complex plane", \emph{Ramanujan J}. \textbf{16}), 105-129, 2008.

\bibitem{Hal1} R.G. Halburd and R.J. Korhonen, ``Nevanlinna theory for the difference operator", \emph{Ann. Acad. Sci. Fenn. Math}.
\textbf{31}, 463-487, 2006.

\bibitem{HKT}R. G. Halburd, R. J. Korhonen and K. Tohge,``Holomorphic curves with shift-invariant hyperplane preimages",
\emph{Trans. Amer. Math. Soc.}, \textbf{366}, 4267--4298, 2014.

\bibitem{HL}Q. Han and F. L\"{u}, ``On the equation $f^n(z)+g^n(z)=e^{\alpha z+\beta}$", \emph{J. Contemp. Math. Anal.,} \textbf{54}(2), 98-102, 2019.

\bibitem{HZ}Z. B. Huang and R. R. Zhang, ``Uniqueness of the differences of meromorphic functions", {\it Analysis Math.}{\bf44}, 461-473, 2018.

\bibitem{Hille}E. Hille, \emph{Ordinary differential equations in the complex domain} (Dover Publication, New York, 1997).

\bibitem{LY}I. Laine and C. C. Yang, ``Clunie theorems for difference and $q$-difference polynomials", \emph{J. London. Math. Soc.}, \textbf{76}, 556--566, 2007.

\bibitem{Latreuch}Z. Latreuch, ``On the existence of entire solutions of certain class of non-linear difference equations", \emph{Mediterr. J. Math.}, \textbf{14}, 1--16, 2017.

\bibitem{L}L. W. Liao and Z. Ye, ``On solutions to nonhomogeneous algebraic differential equations and their application", {\it J. Aus. Math. Soc.}, {\bf97}, 391--403, 2014.

\bibitem{LS}K. Liu and C. J. Song, ``Meromorphic solutions of complex differential-difference equations", \emph{Results Math.}, \textbf{72}, 1759--1771, 2017.

\bibitem{Liu}N.N. Liu, W.R. L\"{u}, T.T. Shen, C.C. Yang, ``Entire solutions of certain type of difference equations", {\it J. Inequal. Appl.,} \textbf{2014}, 63, 2014.

\bibitem{LLX}F. L\"{u}, W. R. L\"{u}, C. P. Li and J. F. Xu, ``Growth and uniqueness related to complex differential and difference equations", \emph{Results Math}. \textbf{74}, 30, 2019.

\bibitem{LvHan}F. L\"{u} and Q. Han, ``On the Fermat-type equation $f^3(z)+ f^3(z +c)=1$",\emph{ Aequat. Math.,} \textbf{91}, 129--136, 2017.

\bibitem{YL}C.C. Yang and I. Laine, ``On analogies between nonlinear difference and differential equations", {\it Proc. Japan Acad. Ser. A.,} {\bf 86}, 10--14, 2010.

\bibitem{Yang}L. Yang, \emph{Value Distribution Theory}(Springer-Verlag, Berlin, 1993)

\bibitem{YY}C.C. Yang and H.X. Yi, \emph{Uniqueness Theory of Meromorphic Functions}(Kluwer, Dordrecht, 2003)

\bibitem{Zhang}F. R. Zhang, N. N. Liu, W. R. L\"{u} and C.C. Yang, ``Entire solutions of certain class of
differential-difference equations", {\it Advances in Difference Equations.} \textbf{2015}, 150, 2015.


\end{thebibliography}
\end{document}